\documentclass[12pt]{article}
\usepackage{amsmath,amssymb}
\usepackage{amsthm}
\usepackage{array}
\newtheorem{thm}{Theorem}[section]
\newtheorem{prop}[thm]{Proposition}
\newtheorem{cor}[thm]{Corollary}
\newtheorem{lem}[thm]{Lemma}

\theoremstyle{note}

\newtheorem{rem}[thm]{Remark}

\def\a{\mathfrak a}

\def\F{\mathfrak F}
\def\M{\mathfrak M}
\def\H{\mathfrak H}

\def\rit#1{{\mbox{\rm #1}}}
\def\modx#1#2{\equiv#1\hspace{-1mm}\mod #2}
\def\nmodx#1#2{\not\equiv#1\hspace{-1mm}\mod #2}
\def\itemx#1{\item[{\rm(#1)}]}
\def\br#1{\{#1\}}

\begin{document}
\title{Special values of generalized $\bf\lambda$ functions at imaginary quadratic points}
\maketitle
\begin{center}
 N{\sc oburo} I{\sc shii}\end{center}
\section{Introduction}
 For a positive integer $N$, let $\Gamma_1(N)$ be the subgroup of $\rit{SL}_2(\mathbf Z)$ defined by
\[
\Gamma_1(N)=\left\{\left. \begin{pmatrix} a & b \\ c & d \end{pmatrix}\in \rit{SL}_2(\mathbf Z)~\right |~ a-1\equiv c \equiv 0 \mod N \right\}.
\]
We denote by $A_1(N)$ the modular function field with respect to $\Gamma_1(N)$. For a positive integer $N\geq 6$, let $\a=[a_1,a_2,a_3]$ be a triple of integers with the properties $0<a_i\leq N/2$ and $a_i\ne a_j$ for any $i,j$. For an element $\tau$ of the complex upper half plane $\mathfrak H$, we denote by $L_\tau$ the lattice of $\mathbf C$ generated by $1$ and $\tau$ and by $\wp(z;L_\tau)$ the Weierstrass $\wp$-function relative to the lattice $L_\tau$. In \cite{II,IK}, we defined a modular function $W_{\a}(\tau)$ with respect to $\Gamma_1(N)$ by
\begin{equation}\label{eqw}
W_{\a}(\tau)=\frac{\wp (a_1/N;\tau)-\wp (a_3/N;\tau)}{\wp (a_2/N;\tau)-\wp (a_3/N;\tau)}.
\end{equation}
This function is one of generalized $\lambda$ functions defined by S.Lang in Chapter 18, \S6 of \cite{LA}. He describes that it is interesting to investigate special values of generalized $\lambda$ functions at imaginary quadratic points and to see if they generate the ray class fields. Here a point of $\H$ is called an imaginary quadratic point if it generates an imaginary quadratic field over $\mathbf Q$. Let $j$ be the modular invariant function. In Theorem 3.6 of \cite{IK}, we showed that $W_\a$ is integral over $\mathbf Z[j]$ under a condition that $a_1a_2a_3(a_1\pm a_3)(a_2\pm a_3)$ is prime to $N$. From this we obtained that values of $W_\a$ at imaginary quadratic points are algebraic integers. Further, we showed in Theorem 5 of \cite{II} that each of the functions $W_{[3,2,1]},W_{[5,2,1]}$ generates $A_1(N)$ over $\mathbf C(j)$. In this article, we study the functions $W_\a$ in the particular case: $a_3=1$. To simplify the notation, henceforth we denote by $\Lambda_{k,\ell}$ the function $W_{[k,\ell,1]}$. We prove that for integers $k,\ell$ such that $1<\ell\ne k<N/2$, the function $\Lambda_{k,\ell}$ generates $A_1(N)$ over $\mathbf C(j)$. This result implies that for an imaginary quadratic point $\alpha$ such that $\mathbf Z[\alpha]$ is the maximal order of the field $K=\mathbf Q(\alpha)$, the value $\Lambda_{k,\ell}(\alpha)$ and $\displaystyle e^{2\pi i/N}$ generate the ray class field of $K$ modulo $N$ over the Hilbert class field $K(j(\alpha))$ of $K$. On the assumption $k(\ell\pm 1)$ is prime to $N$, we prove that $\Lambda_{k,\ell}$ is integral over $\mathbf Z[j]$. Further in the case $\ell =2$,we can weaken the assumption. Let $\delta=(k,N)$ be the greatest common divisor of $k$ and $N$. If we assume either (i)~$\delta=1$ or (ii) $\delta>1,(\delta,3)=1$ and $N/\delta$ is not a power of a prime number, then $\Lambda_{k,2}$ is integral over $\mathbf Z[j]$. In particular the values of $\Lambda_{k,\ell}$ at imaginary quadratic points are algebraic integers. Our results can be extended easily to the functions $W_\a$ in the case that $a_3$ is prime to $N$. See Corollaries \ref{cora},\ref{corb}. Throughout this article, we use the following notation:\newline 
For a function $f(\tau)$ and $A=\begin{pmatrix}a&b\\c&d\end{pmatrix}\in\rit{SL}_2(\mathbf Z)$, the symbols $f[A]_2,f\circ A$ are defined by
\[
f[A]_2=f\left(\frac{a\tau+b}{c\tau+d}\right)(c\tau+d)^{-2},~f\circ A=f\left(\frac{a\tau+b}{c\tau+d}\right).
\]
The greatest common divisor of $a,b\in\mathbf Z$ is denoted by $(a,b)$.
 For an integral domain $R$, $R((q))$ represents the power series ring of a variable $q$ over $R$ and $R[[q]]$ is a subring of $R((q))$ of power series of non-negative order. 

\section{Auxiliary results}
Let $N$ be a positive integer greater than $6$. Put $q=\rit{exp}(2\pi i\tau/N),\zeta=\exp(2\pi i/N)$. For an integer $x$, let
$\{x\}$ and $\mu (x)$ be the integers determined by the following conditions:
\[
\begin{split}
&0\le \{x\}\le \frac N2,\quad \mu (x)=\pm 1,\\
&\begin{cases}\mu(x)=1\qquad &\text{if } x\modx {0,N/2}N,\\
             x\equiv \mu (x)\{x\} \mod N\qquad&\text{otherwise.}
\end{cases}
\end{split}
\]

For an integer $s$ not congruent to $0 \mod N$, let 
\[\phi_s(\tau)=\frac 1{(2\pi i)^2}\wp \left(\frac s N;L_\tau\right)-1/12.
\]
Obviously we have $W_\a=(\phi_{a_1}-\phi_{a_3})/(\phi_{a_2}-\phi_{a_3})$.
 Let $\displaystyle A=\begin{pmatrix}a&b\\c&d\end{pmatrix}\in\rit{SL}_2(\mathbf Z)$. Put $s^*=\mu (sc)sd,u_s=\zeta^{s^*}q^{\{sc\}}$. Then by Lemma 1 of \cite{II}, we have
{\small
\begin{equation}\label{eq1}
\phi_s[A]_2=
\begin{cases}\displaystyle
\frac{\zeta^{s^*}}{(1-\zeta^{s^*})^2}-\sum_{m=1}^{\infty}\sum_{n=1}^{\infty}n(1-\zeta^{s^*n})(1-\zeta^{-s^*n})q^{mnN}&\text{if }\{sc\}=0,\\
\displaystyle\sum_{n=1}^{\infty}n u_s^n-\displaystyle\sum_{m=1}^{\infty}\sum_{n=1}^{\infty}n(1-u_s^n)(1-u_s^{-n})q^{mnN}&\text{otherwise}.
\end{cases}
\end{equation}
}
Therefore, $\phi_s[A]_2\in\mathbf Q(\zeta)[[q]]$ and its order is $\br{sc}$.
If $c\modx 0N$, then by \eqref{eq1} or the transformation formula of $\wp ((r\tau+s)/N;L_\tau)$ in \S2 of \cite{II}, we have
\begin{equation}\label{eqt}
\phi_s[A]_2=\phi_{\br{sd}}.
\end{equation}
The next lemmas and propositions are required in the following sections.
\begin{lem}\label{lem1}
Let $r,s,c,d$ be integers such that $0<r\ne s\leq N/2,~(c,d)=1$. Assume that $\{rc\}=\{sc\}$. Put $r^*=\mu(rc)rd, s^*=\mu(sc)sd$. Then we have $\zeta^{r^*-s^*}\ne 1$. Further if $\{rc\}=\{sc\}=0,N/2$, then $\zeta^{r^*+s^*}\ne 1$.
\end{lem}
\begin{proof}
The assumption $\{rc\}=\{sc\}$ implies that $(\mu(rc)r-\mu(sc)s)c\modx 0N$. If $\zeta^{r^*-s^*}=1$, then $(\mu(rc)r-\mu(sc)s)d\modx 0N$. From $(c,d)=1$, we obtain $\mu(rc)r-\mu(sc)s\modx 0N$. This is impossible, because of $0<r\ne s\leq N/2$. Suppose $\{rc\}=\{sc\}=0,N/2$ and $\zeta^{r^*+s^*}=1$. Then we have $(r+s)c\modx 0N,~(r+s)d\modx 0N$. Therefore $r+s\modx 0N$. This is also impossible.
\end{proof}
\begin{lem}\label{lem2} Let $k\in\mathbf Z$. Put $\delta=(k,N)$. 
\begin{enumerate}
\itemx i For $\ell\in\mathbf Z$, if $\ell$ is divisible by $\delta$, then $(1-\zeta^\ell)/(1-\zeta^k)\in\mathbf Z[\zeta]$.
\itemx {ii} If $N/\delta$ is not a power of a prime number, then $1-\zeta^k$ is a unit of $\mathbf Z[\zeta]$.
\end{enumerate}
\end{lem}
\begin{proof} If $\ell$ is divisible by $\delta$, then there exist an integer $m$ such that $\ell\modx{mk}{N}$. 
Therefore $\zeta^\ell=\zeta^{mk}$ and $(1-\zeta^\ell)$ is divisible by $(1-\zeta^k)$. This shows (i). Let $p_i~(i=1,2)$ be distinct prime factors of $N/\delta$. 
Since $N/p_i=\delta (N/(\delta p_i))$, $1-\zeta^{N/p_i}$ is divisible by $1-\zeta^\delta$. Therefore $p_i~(i=1,2)$ are divisible by $1-\zeta^\delta$. This implies that $1-\zeta^\delta$ is a unit. Because of $(k/\delta,N/\delta)=1$, $1-\zeta^k$ is also a unit .
\end{proof}
From \eqref{eq1} and Lemma~\ref{lem1}, we immediately obtain the following propositions.
\begin{prop}\label{prop1} Let $r,s\in\mathbf Z$ such that $0<r\ne s \leq N/2$.
\begin{enumerate}
\itemx i If $\{rc\},\{sc\}\ne 0$, then 
\[
(\phi_r-\phi_s)[A]_2\equiv \sum_{n=1}^\infty n(u_r^n-u_s^n)+u_r^{-1}q^N-u_s^{-1}q^N \mod q^N\mathbf Z[\zeta][[q]].
\]
\itemx{ii} If $\{rc\}=0$ and $\{sc\}\ne 0$, then
\[(\phi_r-\phi_s)[A]_2\equiv \frac{\zeta^{rd}}{(1-\zeta^{rd})^2}-\sum_{n=1}^\infty nu_s^n-u_s^{-1}q^N \mod q^N\mathbf Z[\zeta][[q]].\]
\itemx {iii} If $\{rc\}=\{sc\}=0$, then 
\[
(\phi_r-\phi_s)[A]_2\equiv \frac{-\zeta^{sd}(1-\zeta^{(r-s)d})(1-\zeta^{(r+s)d})}{(1-\zeta^{rd})^2(1-\zeta^{sd})^2}~\mod q^N\mathbf Z[\zeta][[q]],
\]
\end{enumerate}
\end{prop}
\begin{prop}\label{prop2} Let $r,s\in\mathbf Z$ such that $0<r\ne s \leq N/2$. Put $\ell=\min(\{rc\},\{sc\})$. Then
\[(\phi_r-\phi_s)[A]_2=\theta_{r,s}(A)q^\ell(1+qh(q)),\]
where $h(q)\in\mathbf Z[\zeta][[q]]$ and $\theta_{r,s}(A)$ is a non-zero element of $\mathbf Q(\zeta)$ defined as follows.
In the case $\{rc\}=\{sc\}$,
\[
\theta_{r,s}(A)=\begin{cases}-\zeta^{s^*}(1-\zeta^{r^*-s^*})\quad&\text{if }\ell\ne 0,N/2,\\
           -\zeta^{s^*}(1-\zeta^{r^*-s^*})(1-\zeta^{r^*+s^*})\quad&\text{if }\ell=N/2,\\
\displaystyle\frac{-\zeta^{s^*}(1-\zeta^{r^*-s^*})(1-\zeta^{r^*+s^*})}{(1-\zeta^{r^*})^2(1-\zeta^{s^*})^2}\quad&\text{if }\ell=0.
\end{cases}
\]
In the case $\{rc\}\ne\{sc\}$,assuming that $\{rc\}<\{sc\}$,
\[
\theta_{r,s}(A)=\begin{cases}\displaystyle \zeta^{r^*}\quad&\text{if }\ell\ne 0,\\
\displaystyle\frac{\zeta^{r^*}}{(1-\zeta^{r^*})^2}\quad&\text{if }\ell=0.
\end{cases}
\]
\end{prop}
\section{Generators of $A_1(N)$}
Let $A(N)$ be the modular function field of the principal congruence subgroup $\Gamma (N)$ of level $N$. For a subfield $\F$ of $A(N)$, let us denote by $\F_{\mathbf Q(\zeta)}$ a subfield of $\F$ consisted of all modular functions having Fourier coefficients in $\mathbf Q(\zeta)$.

\begin{thm}\label{th1} Let $k,\ell$ be integers such that $1<\ell\ne k<N/2$. Then we have $A_1(N)_{\mathbf Q(\zeta)}=\mathbf Q(\zeta)(\Lambda_{k,\ell},j)$ 
\end{thm}
\begin{proof} 
 Theorem 3 of Chapter  6 of \cite{LA} shows that $A(N)_{\mathbf Q(\zeta)}$ is a Galois extension over $\mathbf Q(\zeta)(j)$ with Galois group $\rit{SL}_2(\mathbf Z)/\Gamma(N)\{\pm E_2\}$ and $A_1(N)_{\mathbf Q(\zeta)}$ is the fixed field of the subgroup $\Gamma_1(N)\{\pm E_2\}$. Since $\Lambda_{k,\ell}\in A_1(N)_{\mathbf Q(\zeta)}$, to prove the assertion, we have only to show $A\in\Gamma_1(N)\{\pm E_2\}$, for $A\in\rit{SL}_2(\mathbf Z)$ such that $\Lambda_{k,\ell}\circ A=\Lambda_{k,\ell}$. Let $A=\begin{pmatrix}a&b\\ c&d\end{pmatrix}\in\rit{SL}_2(\mathbf Z)$ such that $\Lambda_{k,\ell}\circ A=\Lambda_{k,\ell}$.  Since order of $q$-expansion of $\Lambda_{k,\ell}$ is $0$ and that of $\Lambda_{k,\ell}\circ A$ is $\min(\br{kc},\br{c})-\min(\br{\ell c},\br{c})$ by Proposition \ref{prop2}, we have 
\begin{equation}\label{eq2}
\min(\br{kc},\br{c})=\min(\br{\ell c},\br{c}). 
\end{equation}
From Propositions \ref{prop1} and \ref{prop2}, by considering  power series modulo $q^N$, thus modulo $q^N\mathbf Q(\zeta)[[q]]$, we obtain 
\[
\theta_\ell(\phi_k-\phi_1)[A]_2\equiv \theta_k(\phi_\ell-\phi_1)[A]_2\quad \mod q^N,
\]
where $\theta_\ell=\theta_{\ell,1}(E_2),\theta_k=\theta_{k,1}(E_2)$.
Thus
\begin{equation}\label{eq3}
(\theta_\ell \phi_k-\theta_k\phi_\ell)[A]_2\equiv (\theta_\ell-\theta_k)\phi_1[A]_2\quad \mod q^N.
\end{equation}
If $\theta_\ell=\theta_k$, then we have $(1-\zeta^{k+\ell})(1-\zeta^{k-\ell})=0$. Therefore we have $\theta_\ell\ne\theta_k$.
For an integer $i$, put $u_i=\zeta^{\mu(ic)id}q^{\br{ic}},\omega_i=\zeta^{(\mu(ic)i-\mu(c))d}$. We shall show that $c\modx 0N$. Let us assume $c\nmodx 0N$. Then first of all  we shall prove $\br{\ell c}=\br {kc}=\br c$ by contradiction. Since we can interchange the roles of $\ell,k$, we have only to consider the case $\br{\ell c}\ne\br c$.\vspace{3mm}\newline
(i)~Suppose that $\br{\ell c}<\br c$. Then $\br{kc}=\br{\ell c}, u_k=\omega u_\ell~(\omega=\omega_k/\omega_\ell)$. 
We note $\omega\ne 1$ by Lemma \ref{lem1}. Since order of the power series on the right hand side of \eqref{eq3} is $\br c\ne 0$, we see $\br{\ell c}\ne 0$, because if $\br{\ell c}= 0$, then order of the series on the left hand side is a multiple of $N$. Further the coefficient $\theta_\ell\omega-\theta_k$ of $u_\ell$ of the series on the left hand side should be $0$. Thus from Proposition \ref{prop1} and \eqref{eq3}, we obtain
\[(\phi_k-\omega \phi_\ell)[A]_2\equiv (1-\omega)\phi_1[A]_2\quad\mod q^N.\]
This gives
\[
\sum_{n\geq 2}n\omega(\frac{\omega^{n-1}-1}{\omega-1})u_\ell^n-\omega^{-1}(\omega+1)u_\ell^{-1}q^N\equiv -\sum_{n\geq 1}nu_1^n-u_1^{-1}q^N~\mod q^N.
\]
 Since $\br c\leq N-\br c<N-\br{\ell c}$,we have $2\br{\ell c}=\br c$. If $\br c=N/2$, then $\br{\ell c}=0,N/2$ according to the parity of $\ell$. Therefore we know $\br c\ne N/2$ and by comparing the coefficients of $q^{\br c}$ on both sides, we have $2\omega\zeta^{2\mu(\ell c)\ell d}=-\zeta^{\mu(c)c}$. This gives a contradiction.\vspace{3mm}\newline
(ii) Suppose that $\br{\ell c}>\br c$. The congruence \eqref{eq3} implies $\br{kc}=\br c$ and 
\[
\theta_k\phi_\ell[A]_2\equiv(\theta_\ell\phi_k-(\theta_\ell-\theta_k)\phi_1)[A]_2\quad\mod q^N.
\]
Since $\br c\ne 0,N/2$ and $u_k=\omega_ku_1$, we have $\theta_k=(1-\omega_k)\theta_\ell$ and, noting $\omega_k\ne 1$,
\[
\sum_{n\geq 1}nu_\ell^n+u_\ell^{-1}q^N\equiv\sum_{n\geq 2}n\left(1-\frac{\omega_k^n-1}{\omega_k-1}\right)u_1^n+(\omega_k^{-1}+1)u_1^{-1}q^N
\quad\mod q^N.
\]
Since $N-\br c>N-\br{\ell c}\geq\br{\ell c}$, we have $2\br c=\br{\ell c}$. By comparing the coefficients of first terms on both sides, we obtain $2\omega_k\zeta^{2\mu(c)d}=-(\zeta^{\ell d}+\zeta^{- \ell d})$ in the case $\br{\ell c}=N/2$ and  $2\omega_k\zeta^{2\mu(c)d}=-\zeta^{\ell d}$ in the case $\br{\ell c}\ne N/2$. 
In the former case,  we have $|\cos 2\pi \ell d/N|=1$. Therefore $\ell c\modx 0{N/2}$ and $\ell d\modx 0{N/2}$. 
Since $(c,d)=1$, we know $\ell\modx 0{N/2}$. This is impossible. In the latter case, clearly we have a contradiction. Therefore we have $\br{\ell c}=\br {kc}=\br c$. Assume that $\br{\ell c}=\br {kc}=\br c$. Then $u_k=\omega_ku_1,u_\ell=\omega_\ell u_1$. By congruence \eqref{eq3}, by comparing the coefficients of $u_1$ on both sides, we have $\theta_\ell(\omega_k-1)=\theta_k(\omega_\ell-1)$. Since $\omega_k,\omega_\ell,\omega_k/\omega_\ell \ne 1$, from \eqref{eq3}, we obtain
\[
2u_1^2+(\omega_k\omega_\ell)^{-1}u_1^{-1}q^N\equiv  0\quad \mod u_1^3.
\]
This gives a contradiction. Hence we have $c\modx 0N$. Hereafter we assume that $k>\ell$, if necessary by considering $1/\Lambda_{k,\ell}=\Lambda_{\ell,k}$ instead of $\Lambda_{k,\ell}$.
Let $c\modx 0N$. Then $(d,N)=1$. By \eqref{eqt}, we have $\Lambda_{k,\ell}\circ A=\frac{\phi_{\br{kd}}-\phi_{\br{d}}}{\phi_{\br{\ell d}}-\phi_{\br{d}}}$. From now on, to save labor, we put $r=\br{\ell d},s=\br{kd},t=\br d$. Since $r,s,t$ are distinct from each other and $\min(s,t)=\min(r,t)$, we have $r,s,t\ne 0,N/2$ and $t<r,s$. We have only to prove $t=1$. Let us assume $t>1$.
Let $T=\begin{pmatrix}1&0\\1&1\end{pmatrix}$. Then 
\begin{equation}\label{eq5}
\Lambda_{k,\ell}\circ T=\Lambda_{k,\ell}\circ AT=\left(\frac{\phi_s-\phi_t}{\phi_r-\phi_t}\right)\circ T.
\end{equation}
If $i$ is an integer such that $0<i<N/2$, then $\mu(i)=1,\br i=i$. Let $u=\zeta q$. Then  
\begin{equation}\label{eq6}
\phi_i[T]_2\equiv \sum_n nu^{i n}+u^{N-i} \mod q^N.
\end{equation}
From \eqref{eq5},
\[(\phi_r\phi_1+\phi_s\phi_\ell+\phi_t\phi_k)[T]_2=(\phi_t\phi_\ell+\phi_s\phi_1+\phi_r\phi_k)[T]_2.
\]
Since $s+\ell,t+k,r+k>t+\ell$ and order of $\phi_i[T]_2$ is $i$, we have
\[\phi_r\phi_1[T]_2-\phi_s\phi_1[T]_2\equiv \phi_t\phi_\ell[T]_2~\mod u^{t+\ell+1}.\]If $s<r$, then $s+1=t+\ell$. However in this case the coefficients of $u^{t+\ell}$ on both sides are distinct. Therefore $r<s,~r+1=t+\ell$. Since $t+\ell\leq s<N/2$, we know that $N>2t+2\ell$. 
By \eqref{eq6} and by the inequality relations that $t\geq 2,~k>\ell$,$r=t+\ell -1,s\geq t+\ell,N>2t+2\ell$, we have modulo $u^{t+\ell+2}$,
\[ 
\begin{split}
&(\phi_r\phi_1)[T]_2\modx{u^{t+\ell}+2u^{t+\ell+1}}{u^{t+\ell+2}},~(\phi_s\phi_\ell)[T]_2\modx{0}{u^{t+\ell+2}},\\
&(\phi_t\phi_k)[T]_2\modx{u^{t+k}}{u^{t+\ell+2}},~(\phi_t\phi_\ell)[T]_2\modx{u^{t+\ell}}{u^{t+\ell+2}},\\
&(\phi_s\phi_1)[T]_2\modx{u^{s+1}}{u^{t+\ell+2}},~(\phi_r\phi_k)[T]_2\modx{0}{u^{t+\ell+2}}.
\end{split}
\]
Therefore we obtain a congruence:
\[2u^{t+\ell +1}+u^{t+k}\modx{u^{s+1}}{u^{t+\ell +2}}.\]
The coefficients of $u^{t+\ell +1}$ on both sides cannot be equal. Hence we have a contradiction.
\end{proof}
\begin{cor}\label{cora}
 Let $W_\a$ be the function defined by \eqref{eqw}.
If $a_1,a_2\ne N/2$ and $(a_3,N)=1$, then $A_1(N)_{\mathbf Q(\zeta)}=\mathbf Q(\zeta)(W_\a,j)$
\end{cor}
\begin{proof}
Let $M \in\rit{SL}_2(\mathbf Z)$ such that $M\equiv\begin{pmatrix}a_3^{-1}&0\\ 0&a_3\end{pmatrix}\mod N$. By \eqref{eqt}, we know that $W_\a=\Lambda_{k,\ell}\circ M$, where $k,\ell\in\mathbf Z$ such that $a_1=\br{ka_3},a_2=\br{\ell k}$ and $1<\ell\ne k<N/2$. Let $A\in \rit{SL}_2(\mathbf Z)$. If $W_\a\circ A=W_\a$, then $\Lambda_{k,\ell}\circ (MAM^{-1})=\Lambda_{k,\ell}$. Therefore we have $MAM^{-1}\in\Gamma_1(N)\{\pm E_2\}$. Since $M$ is a normalizer of $\Gamma_1(N)$, $A\in\Gamma_1(N)\{\pm E_2\}$. This shows our assertion.
\end{proof}
\section{Values of $\Lambda_{k,\ell}$ at imaginary quadratic points}
In this section, we shall study values of $\Lambda_{k,\ell}$ at imaginary quadratic points.
\begin{prop}\label{prop3a} Let $k,\ell$ be integers such that $1<\ell\ne k\leq N/2$.  Assume $(k(\ell\pm 1),N)=1$. Then for $A\in\rit{SL}_2(\mathbf Z)$,we have
\[\Lambda_{k,\ell}\circ A\in\mathbf Z[\zeta]((q)).\]
\end{prop}
\begin{proof} Put $A=\begin{pmatrix}a&b\\c&d\end{pmatrix}$. Proposition \ref{prop2} shows 
 \[\Lambda_{k,\ell}\circ A=\omega f(q),\]
 where $\omega=\theta_{k,1}(A)/\theta_{\ell,1}(A)$ and $f$ is a power series in $\mathbf Z[\zeta]((q))$. Therefore it is sufficient to prove that $\omega\in \mathbf Z[\zeta]$. First we consider the case $\{c\}\ne 0$. By the assumption, we know $\{\ell c\}\ne\{c\}$. Proposition\ref{prop2} implies that $\theta_{\ell,1}(A)^{-1} \in \mathbf Z[\zeta]$. Since $(k,N)=1$, $\br{kc}\ne 0$. Thus $\theta_{k,1}(A)\in\mathbf Z[\zeta]$. Hence we have $\omega\in\mathbf Z[\zeta]$. 
Next consider the case $\{c\}=0$. Then we have $\{c\}=\{\ell c\}=\{kc\}=0,\mu(c)=\mu(\ell c)=\mu(kc)=1$, $(d,N)=1$ and 
\[
\omega=\left(\frac{1-\zeta^{\ell d}}{1-\zeta^{kd}}\right)^2\cdot\frac{(1-\zeta^{(k-1)d})(1-\zeta^{(k+1)d})}{(1-\zeta^{(\ell-1)d})(1-\zeta^{(\ell+1)d})}.
\]
Using the assumption, Lemma \ref{lem2} (i) shows  
\[\displaystyle\frac{1-\zeta^{\ell d}}{1-\zeta^{kd}},\frac{1-\zeta^{(k-1)d}}{1-\zeta^{(\ell-1)d}},\frac{1-\zeta^{(k+1)d}}{1-\zeta^{(\ell+1)d}}\in\mathbf Z[\zeta].\]
 Hence we obtain $\omega\in\mathbf Z[\zeta]$. 
\end{proof}
\begin{prop}\label{prop3} Let $k$ be an integer such that $2<k<N/2$. Put $\delta=(k,N)$. Assume either \rit{(i)} $\delta=1$ or \rit{(ii)} $\delta>1,(\delta,3)=1$ and $N/\delta$ is not a power of a prime number. Then for $A\in\rit{SL}_2(\mathbf Z)$,we have
\[\Lambda_{k,2}\circ A\in\mathbf Z[\zeta]((q)).\]
\end{prop}
\begin{proof} Similarly in the proof of Proposition \ref{prop3a}, we have only to prove  $\omega=\theta_{k,1}(A)/\theta_{2,1}(A)\in \mathbf Z[\zeta]$. First we consider the case $\{c\}\ne 0$. Let $\{2c\}\ne\{c\}$
. By (ii) of Proposition \ref{prop2}, we see $\theta_{2,1}(A)^{-1} \in \mathbf Z[\zeta]$. 
Further if $\{kc\}\ne0$, then $\theta_{k,1}(A)\in\mathbf Z[\zeta]$. If $\{kc\}=0$,then $\delta>1$ and $c\modx0{N/\delta}$. Therefore $\zeta^{kd}$ is a primitive $N/\delta$-th root of unity. The assumption (ii) shows $1-\zeta^{kd}$ is a unit. Thus $\theta_{k,1}(A)\in\mathbf Z[\zeta]$. Hence we have $\omega\in\mathbf Z[\zeta]$. Let $\{2c\}=\{c\}$. 
Then, since $\{c\}\ne 0$, we have $N\modx 03,~(k,3)=1$ and $\{c\}=\{2c\}=\{kc\}=N/3$, $\mu(2c)=-\mu(c)$, $\mu(kc)=(\frac k3)\mu(c)$, where $(\frac *3)$ is the Legendre symbol. 
By the same proposition, we know that $\displaystyle\omega=(1-\zeta^{(\mu(kc)k-\mu(c))d})/(1-\zeta^{-3\mu(c)d})$. Since $\mu(kc)k-\mu(c)\modx 03$, we have $\omega\in\mathbf Z[\zeta]$. 
Next consider the case $\{c\}=0$. Then we have $\{c\}=\{2c\}=\{kc\}=0,\mu(c)=\mu(2c)=\mu(kc)=1$, $(d,N)=1$ and 
\[
\omega=\left(\frac{1-\zeta^{2d}}{1-\zeta^{kd}}\right)^2\cdot\frac{(1-\zeta^{(k-1)d})(1-\zeta^{(k+1)d})}{(1-\zeta^d)(1-\zeta^{3d})}.
\]
If $\delta =1$, then $(kd,N)=1$. If $\delta\ne 1$, then the assumption (ii) implies $(1-\zeta^{kd})$ is a unit. Therefore $(1-\zeta^{2d})/(1-\zeta^{kd})\in\mathbf Z[\zeta]$.  If $N\not\equiv 0\mod 3$, then since $(3d,N)=1$, we know \[\frac{(1-\zeta^{(k-1)d})(1-\zeta^{(k+1)d})}{(1-\zeta^d)(1-\zeta^{3d})}\in\mathbf Z[\zeta].\]
 If $N\modx 03$, then $(k,3)=1$ and one of $k+1,k-1$ is divisible by $3$. Lemma \ref{lem1} (i) gives  
\[\displaystyle\frac{(1-\zeta^{(k-1)d})(1-\zeta^{(k+1)d})}{(1-\zeta^d)(1-\zeta^{3d})}\in\mathbf Z[\zeta].\]
 Hence we obtain $\omega\in\mathbf Z[\zeta]$. 
\end{proof}
To study the modular equation of $\Lambda_{k,\ell}$ over $\mathbf C(j)$, we construct a transversal $R$ of the coset decomposition of $\rit{SL}_2(\mathbf Z)$ by $\Gamma_1(N)\{\pm E_2\}$, where $E_2$ is the unit matrix. For $v \in (\mathbf Z/N\mathbf Z)^\times/\br{\pm 1}$, take $M_v\in\rit{SL}_2(\mathbf Z)$ so that $M_v\equiv \big(\begin{smallmatrix}v^{-1}&0\\0&v\end{smallmatrix}\big)\mod N$. For a positive divisor $t$ of $N$, let $\Theta_t$ be a set of $\varphi((t,N/t))$ integers $u$ such that $u$ is prime to $t$ and runs over a transversal of the factor group $(\mathbf Z/(t,N/t)\mathbf Z)^\times$. For an integer $u$ prime to $t$ and an integer $k$ with the property $ku\modx 1t$, consider a matrix in $\rit{SL}_2(\mathbf Z)$ 
\[
B(t,u,k)=\begin{pmatrix}u&(uk-1)/t\\t&k\end{pmatrix}.
\]
 We denote by $\M_{\Theta_t}$ the set of matrices 
$$\{B(t,u,k)~|~u\in\Theta_t,~k \rit{ mod } N/t,uk\modx 1t\}.$$ 
Lemma 3.1 in \cite{IK} shows a set of matrices 
\[R=\{M_v B~|~v\in(\mathbf Z/N\mathbf Z)^\times/\br{\pm 1},B\in\underset{t|N}{\cup}\M_{\Theta_t}\}
\] is a transversal of the coset decomposition of $\rit{SL}_2(\mathbf Z)$ by $\Gamma_1(N)\{\pm E_2\}$. For an integer $h$ prime to $N$, let $h^*$ be an integer such that $hh^*\modx 1N$. Put , for the set $\Theta_t$, 
\[
h\Theta_t=\{h^*u|u\in\Theta_t\}.
\]
Then obviously a set of matrices $\{M_v B~|~v\in(\mathbf Z/N\mathbf Z)^\times/\br{\pm 1},B\in\underset{t|N}{\cup}\M_{h\Theta_t}\}$ is also a transversal of the coset decomposition. For details, see \S 3 of \cite{IK}. Let $\sigma_h$ be the automorphism of $\mathbf Q(\zeta)$ defined by $\zeta^{\sigma_h}=\zeta^h$. On a power series $f=\sum_ma_mq^m$ with $a_m\in\mathbf Q(\zeta)$, $\sigma_h$ acts by $f^{\sigma_h}=\sum_m a_m^{\sigma_h}q^m$.
From Lemma 3.2 of \cite{IK}, we obtain
\begin{equation}\label{eqa}
(\Lambda_{k,\ell}\circ M_v B(t,u,k))^{\sigma_h}=\Lambda_{k,\ell}\circ M_v B(t,h^*u,hk).
\end{equation} 
\begin{rem} We do not need to separate the case $\br {st}=0$ from the case $\br {st}\ne 0$ in Lemma 3.2 of \cite{IK}. Therefore we can omit the assumption $a_1a_2a_3$ is prime to $N$ in Proposition 3.4,\cite{IK}.
\end{rem}
\begin{thm}\label{th2} Let the assumption be the same as in Proposition \ref{prop3a}. Further in the case $\ell=2$, let the assumption be the same as in Proposition \ref{prop3}. Then $\Lambda_{k,\ell}$ is integral over $\mathbf Z[j]$. 
\end{thm}
\begin{proof}
Let $R$ be the above set. Consider a modular equation $\Phi(X,j)=\prod_{A\in R}(X-\Lambda_{k,\ell}\circ A)$. Since $\Lambda_{k,\ell}\circ A$ has no poles in $\H$ and $\Lambda_{k,\ell}\circ A\in\mathbf Z[\zeta]((q))$ by Propositions \ref{prop3a} and \ref{prop3}, the coefficients of $\Phi(X,j)$ are polynomials of $j$ with coefficients in $\mathbf Z[\zeta]$. By \eqref{eqa}, for every automorphism $\sigma$ of $\mathbf Q(\zeta)$, the correspondence: $\Lambda_{k,\ell}\circ A\rightarrow (\Lambda_{k,\ell}\circ A)^\sigma$ induces a permutation on the set $\{\Lambda_{k,\ell}\circ A~|~A\in R\}$. Therefore $\Phi(X,j)\in\mathbf Z[j]$.
\end{proof}
\begin{thm} Let the assumption be the same as in Theorem \ref{th2}. Let $\alpha$ be an imaginary quadratic point. Then $\Lambda_{k,\ell}(\alpha)$ is an algebraic integer.
\end{thm}
\begin{proof}
 Since $j(\alpha)$ is an algebraic integer (see \cite{C1},Theorem 10.23) and $\Lambda_{k,\ell}(\alpha)$ is integral over $\mathbf Z[j(\alpha)]$, $\Lambda_{k,\ell}(\alpha)$ is an algebraic integer. 
\end{proof}
\begin{cor} Let $A\in\rit{SL}_2(\mathbf Z)$. Let the assumption be the same as in Theorem \ref{th2}. Then the values of the function $\Lambda_{k,\ell}\circ A$ at imaginary quadratic points are algebraic integers. In particular, the function
\[
\frac{\wp (k\tau/N;\tau)-\wp (\tau/N;\tau)}{\wp (\ell\tau/N;\tau)-\wp (\tau/N;\tau)}
\]
takes algebraic, integral values at imaginary quadratic points. 
\end{cor}
\begin{proof}
Let $\alpha$ be an imaginary quadratic point. Then, $A(\alpha)$ is an imaginary quadratic point. Therefore, we have the former part of the assertion. If we put $\displaystyle A=\begin{pmatrix}0&-1\\1&0 \end{pmatrix}$, then from the transformation formula of $\wp ((r\tau+s)/N;L_\tau)$ in \S2 of \cite{II}, we obtain the latter part.
\end{proof} 
\begin{cor}\label{corb} Let $W_\a$ be the function defined by \eqref{eqw}. If $a_1a_3(a_2\pm a_3)$ is prime to $N$, then $W_\a$ is integral over $\mathbf Z[j]$ and values of $W_\a$ at imaginary quadratic points are algebraic integers.
\end{cor}
\begin{proof} 
Let $M \in\rit{SL}_2(\mathbf Z)$ such that $M\equiv\begin{pmatrix}a_3^{-1}&0\\ 0&a_3\end{pmatrix}\mod N$. Then $W_\a=\Lambda_{k,\ell}\circ M$, where $k,\ell\in \mathbf Z$ such that $a_1=\br{ka_3},a_2=\br{\ell a_3}$. The assumption implies $k(\ell\pm 1)$ is prime to $N$. By Theorem \ref{th2}, for $A\in\rit{SL}_2(\mathbf Z)$, $\Lambda_{k,\ell}\circ A$ is integral over $\mathbf Z[j]$. Hence $W_\a$ is integral over $\mathbf Z[j]$.
\end{proof}
We obtain the following theorem from the Gee-Stevenhagen theory in \cite{GA} and \cite{GAS}. See also Chapter 6 of \cite{SG}.
\begin{thm}
Let $k,\ell\in\mathbf Z$ such that $1<\ell\ne k<N/2$. Let $K$ be an imaginary quadratic field with the discriminant $D$. Then the ray class field of $K$ modulo $N$ is generated by $\Lambda_{k,\ell}((D+\sqrt D)/2)$ and $\zeta$ over the Hilbert class field of $K$.
\end{thm}
\begin{proof}
The assertion is deduced from Theorems 1 and 2 of \cite{GA} and Theorem \ref{th1}.
\end{proof}

\end{document}